\documentclass[reqno]{amsart}
\usepackage{amsmath,amssymb,cite}
\usepackage[mathscr]{euscript}
\usepackage{mathtools}
\usepackage{bbm}
\mathtoolsset{showonlyrefs}
\usepackage{xcolor}

\usepackage{ulem}

\definecolor{bblue}{rgb}{.2,0.2,.8}

\theoremstyle{plain}
\newtheorem{theorem}{Theorem}[section]

\theoremstyle{definition}

\theoremstyle{remark}

\newcommand\one{\mathbbm{1}}

\numberwithin{equation}{section}
\numberwithin{theorem}{section}

\def\be{\begin{equation}}
\def\ee{\end{equation}}
\def\bp{\begin{pmatrix}}
\def\ep{\end{pmatrix}}
\def\bea{\begin{eqnarray}}
\def\eea{\end{eqnarray}}

\def\\{\par\medskip}

\let\0=\noindent

\renewcommand{\epsilon}{\varepsilon}

\title[Solvable SNS]{Solvable stationary non equilibrium states}

\author{G. Carinci}
\address{\noindent Gioia Carinci \hfill\break\indent 
	FIM, Universit\`a di Modena e Reggio Emilia
	\hfill\break\indent 
	42121 Modena, Italy
}
\email{gioia.carinci@unimore.it}

\author{C. Franceschini}
\address{\noindent Chiara Franceschini \hfill\break\indent 
	FIM, Universit\`a di Modena e Reggio Emilia
	\hfill\break\indent 
	42121 Modena, Italy
}
\email{chiara.franceschini@unimore.it}

\author{D. Gabrielli}
\address{\noindent Davide Gabrielli \hfill\break\indent 
	DISIM, Universit\`a dell'Aquila
	\hfill\break\indent 
	67100 Coppito, L'Aquila, Italy
}
\email{davide.gabrielli@univaq.it}

\author{C. Giardin\`a}
\address{\noindent Cristian Giardin\`a \hfill\break\indent 
	FIM, Universit\`a di Modena e Reggio Emilia
	\hfill\break\indent 
	42121 Modena, Italy
}
\email{cristian.giadrina@unimore.it}

\author{D. Tsagkarogiannis}
\address{\noindent Dimitrios Tsagkarogiannis\hfill\break\indent 
	DISIM, Universit\`a dell'Aquila
	\hfill\break\indent 
	67100 Coppito, L'Aquila, Italy
}
\email{dimitrios.tsagkarogiannis@univaq.it}

\begin{document}

\begin{abstract}
We consider the one dimensional boundary driven harmonic model and its continuous version, both introduced in \cite{FGK}. By combining duality and integrability the authors of \cite{FG} obtained the invariant measures in a combinatorial representation. Here we give an integral representation of the invariant measures which
turns out to be a convex combination of inhomogeneous product of geometric distributions for the discrete model and a convex combination of inhomogeneous product of exponential distributions for the continuous one. The mean values of the geometric and of the exponential variables are distributed according to the order statistics of i.i.d. uniform random variables on a suitable interval fixed by the boundary sources. The result is obtained solving exactly the stationary condition written in terms of the joint generating function. The method has an interest in itself and can be generalized to study other models. We briefly discuss some applications.
\end{abstract}

\noindent


\maketitle
\thispagestyle{empty}

\section{Introduction}

Stationary non equilibrium states (SNS) have a rich and complex structure. A natural way to generate a SNS using stochastic interacting particle systems is to put the system, that is evolving on a lattice, in contact with external sources. 
This is a toy model for a thermodynamic system with external reservoirs. The Markov process obtained with this procedure is typically non-reversible when the reservoirs have different parameters
and its invariant measure is the SNS. Due to the non reversibility, such measure is typically difficult to be computed and has long range correlations \cite{MFT,Derr}.

From a {\it macroscopic} point of view, for a few one dimensional solvable models it is possible to get a description of the fluctuations of the SNS by an exact computation of the density large deviations rate functional. This is obtained either by using combinatorial representation of the invariant measure \cite{Derr} or by the variational dynamic approach of Macroscopic Fluctuation Theory (MFT) \cite{MFT}.
Among the solvable models there are the symmetric exclusion process (SEP) and the Kipnis-Marchioro-Presutti (KMP) model \cite{KMP,BGL} and more generally all the models having a constant diffusion and a quadratic mobility in the hydrodynamic scaling limit.
Due to the presence of long range correlations, the rate functionals are non-local and can be written in terms of the maximization  (for SEP) or minimization  (for KMP) of an auxiliary function. 
A problem of interest is the interpretation of the auxiliary function. 
In the case of the KMP model it has been conjectured in \cite{BDGJL-int} that the auxiliary function can be interpreted as a hidden temperature and the minimization as a contraction principle. This conjecture is solved in \cite{DFG} where a joint energy-temperature dynamics has been constructed; as a consequence the invariant measure of the boundary driven case is written as a convex combination of inhomogeneous product of exponential distributions whose mean values are distributed according to the invariant measure of an auxiliary opinion model.

From a {\it microscopic} point of view, before the most recent developments, to our knowledge, there were essentially a few models with long-range correlations for which the description of the stationary measure was explicit. This is the class of open
exclusion type processes, for which it is available a {matrix product ansatz} (see  \cite{DEHP, SS}). 
It was exactly this explicit knowledge that made possible to obtain the density large
deviation function by a microscopic computation \cite{DLS} and 
{then to verify the agreement with the variational structure of MFT \cite{BDGJL1}}.  Since MFT is believed to have a large degree of universality,
as the theory describing fluctuations in diffusive systems, it is therefore important to have additional models of which the SNS is known.
Furthermore, for stationary non-equilibrium states a general structure does not exist as it is the case for equilibrium, where one has instead the Boltzmann-Gibbs distribution.

In a series of recent works \cite{FGK,FFG,FG}, two new  integrable models have been introduced. 
These are the family of  harmonic models, a class of interacting particle systems,  and a suitable continuous version, that can be interpreted as a model for heat conduction. The latter is obtained as a scaling limit of the discrete one. Both families of models are parametrized by a (spin) value $s>0$.
 The integrability of the systems relies on an algebraic description of the generator and the link with integrable systems in quantum spin chains, as is the case also for the class of exclusion processes (see e.g. \cite{SS}). Besides sharing the same algebraic description, these two models are also in a duality relation via a moment duality function \cite{FFG}.
Both models are of zero range type, i.e., the rate at which particles or energy is transferred from one site to another depends just on the number of particles or amount of energy present on the departure site. 
However, differently from the classic zero range models, here there are transitions of multiple particles and the boundary driven SNS are not of product type. 
In Corollary 2.9 of \cite{FG} a closed formula\footnote
{When constructing the mixed measures of this paper, we used this formula to check that a product of geometric with
 mixing measure given by the ordered statistics of i.i.d. uniforms was indeed reproducing the correlation functions in \cite{FG}.}

of combinatorial nature for the stationary measure has been obtained for the class of harmonic models.
The derivation relies on  techniques inspired by integrable systems  and is based on  a direct mapping between non-equilibrium and equilibrium \cite{TKL,FGK2}.
A similar study has been done in  \cite{FFG} for the family of integrable heat conduction models, for which moments of the stationary measure have been found via stochastic duality. For both classes of models an explicit description of the long-range correlations has been shown. As shown in \cite{CGT} all these models have constant diffusion and quadratic mobility and are therefore good candidates for having a mixture of product distribution as invariant measure.

In this paper we provide the probabilistic description of the SNS of these two models from a microscopic perspective.
We consider  the special case  $s=1/2$ for the (spin) value.  We prove that for this pair of models the invariant measure can be written, like for the KMP model, as a mixture of  products of inhomogeneous distributions. Furthermore, for the models considered here, the mixing measure can be explicitly characterized in terms of the order statistics of i.i.d. uniform random variables.  This probabilistic interpretation  sheds light on how the structure of long-range correlations  of the SNS is  rooted in the correlated structure of the mixing measure.

In the harmonic model of parameter $s=1/2$ considered here, at each site of a graph there is a non-negative integer number of particles. When on a vertex $x$ there are $\eta_x$ particles, then $k\leq \eta_x$ particles can jump across each edge exiting from $x$ with rate $1/k$. We consider a one-dimensional lattice with left and right extrema coupled to reservoirs having densities $0<\rho_A\leq \rho_B<+\infty$. When $\rho_A = \rho_B$ the model is reversible and its invariant measure is of product type with each marginal being geometric with mean equal to the density of the external reservoirs. When $\rho_A < \rho_B$ we prove that the invariant measure is a mixture of inhomogeneous product of geometric distributions.
The law of the mean values of the inhomogeneous geometric distributions is the order statistics of independent uniform random variables in the interval $[\rho_A , \rho_B]$. This is  a natural representation, 
since the computation of the integral over the hidden parameters does not give a transparent expression, being written in terms of hypergeometric functions. For the continuous model we have a similar representation, the heat baths attached at the end points of the bulk have temperatures $0<T_A\leq T_B<+\infty$ and the geometric distributions have to be substituted by the exponential ones.

Concerning the methodology, our result is proved writing the stationarity condition of the master equation 
in terms of the joint generating function. This allows a direct verification of the mixed structure of the SNS via a telescopic property.
In this paper we apply the method just to two models in order to give a direct and clear presentation. 
We plan to give a systematic study in the future.
We believe the mixed structure with random temperatures/chemical potentials of the stationary measure to be common to several open models of interacting particles, like for example the exclusion process. This is related also to the fact that the corresponding large deviations rate functionals can be written equivalently in terms of infimum (see \cite{BDGJL-int, DEL, EnauD}). 
See \cite{FGC} for results in this direction for the symmetric exclusion process.

\bigskip

{\bf Note added:} 
After this article was submitted reference \cite{CFFGR} appeared on the arXiv. 
It contains the mixed measure for the harmonic model with general spin s, which
is obtained by a constructive approach that allows to identify the ordered Dirichlet process
as the mixing measure.
It further contains a direct proof that the measure for s=1/2 derived here coincides with the one derived in \cite{FG}
(see Appendix A of \cite{CFFGR} for the comparison).

\section{The discrete Harmonic model with parameter $s=1/2$} \label{sec2}

\subsection{The model} \label{sec11}
We consider a one-dimensional lattice consisting of $N$ sites (the bulk) $\Lambda_N:=\{1,\dots, N\}$ and 
two ghost lattice sites (the boundaries) $\partial \Lambda_N:=\{0, N+1\}$ to which we associate two  parameters $0<\beta_A <\beta_B<1$, respectively. On each lattice site we can have an arbitrarily large number of particles and we denote by $\eta_x\in \mathbb N_0$ the number (possibly zero) of particles at $x\in \Lambda_N$. We consider a continuous-time Markov chain $\{\eta(t), \ t\ge 0\}$ whose state space is the set  
$\Omega_N=\mathbb N_0^{\Lambda_N}$ of configurations $\eta=(\eta_1,\dots ,\eta_N)$, with  $\eta_x\in \mathbb N_0$ being the number of particles at site $x\in \Lambda_N$. The stochastic dynamics has a bulk and a boundary part which are described in terms of the generator $L_N$ defined below. For any $x\in \Lambda_N$, we denote by $\delta_x \in \Omega_N$ the configuration defined by $\delta_x(y)=0$ when $y\neq x$ and $\delta_x(x)=1$.
We have 
\begin{equation}\label{generator}
L_N:=L_N^{\textrm{bulk}}+L_N^{\textrm{bound}}.
\end{equation}
The bulk generator applied to bounded functions reads:
\begin{equation}
\label{genbulk}
L_N^{\textrm{bulk}}f(\eta)= \sum_{\substack{x,y\in \Lambda_N \\ |x-y|=1}} \sum_{k=1}^{\eta_x}\frac 1k \left[f(\eta-k\delta_x+k\delta_y)-f(\eta)\right]\,.
\end{equation} 
Furthermore, the boundary part which encodes the interaction with the reservoirs
is given by:
\begin{align}
\label{genbound}
L_N^{\textrm{bound}}f(\eta) & = \sum_{k=1}^{\eta_1}\frac 1k \left[f(\eta-k\delta_1)-f(\eta)\right]+
\sum_{k=1}^{\infty}\frac{\beta_A^k}{k}\left[f(\eta+k\delta_1)-f(\eta)\right]\nonumber \\
& + \sum_{k=1}^{\eta_N}\frac 1k \left[f(\eta-k\delta_{N})-f(\eta)\right]+
\sum_{k=1}^{\infty}\frac{\beta_B^k}{k}\left[f(\eta+k\delta_{N})-f(\eta)\right]\,.
\end{align}

\subsection{Invariant measure}

For a generic measure $\mu$ on $\Omega_N$
the stationarity condition $\mu L_N=0$
reads as follows:

\begin{align}\label{stazcondN}
&\mu(\eta)\left[\sum_{k=1}^{\infty}\frac{\beta_{A}^k}{k}+\sum_{x=1}^{N}\sum_{k=1}^{\eta_x}\frac 2k+\sum_{k=1}^{\infty}\frac{\beta_{B}^k}{k}\right]\\
&
=\sum_{k=1}^{\eta_1}\mu(\eta-k\delta_1)\frac{\beta_{A}^k}{k} + \sum_{k=1}^{\infty}\mu(\eta+k\delta_1)\frac 1k \nonumber \\
&
+\sum_{x=1}^{N-1}\sum_{k=1}^{\eta_{x+1}}\mu(\eta+k\delta_x-k\delta_{x+1})\frac 1k
 +\sum_{x=2}^{N}\sum_{k=1}^{\eta_{x-1}}\mu(\eta+k\delta_x-k\delta_{x-1})\frac 1k \nonumber \\
&
+\sum_{k=1}^{\eta_N}\mu(\eta-k\delta_{N})\frac{\beta_{B}^k}{k}+\sum_{k=1}^{\infty}\mu(\eta+k\delta_{N})\frac 1k\,.
\end{align}
Let $\mathcal G_m(k)=\frac{1}{1+m}\left(\frac{m}{1+m}\right)^k$, $k=0,1,\dots$, be a geometric distribution of mean $m$. Given $\underline m=(m_1,\dots ,m_{N})$ and $\underline k =(k_1, \dots ,k_{N})$ we denote by
$\mathcal G_{\underline m}(\underline k):=\prod_{x=1}^{N}\mathcal G_{m_x}(k_x)$.
Given $0<\beta_A<\beta_B <1$ we call $\rho_A:=\frac{\beta_A}{1-\beta_A}$, $\rho_B:=\frac{\beta_B}{1-\beta_B}$ and introduce
$O^{\rho_A,\rho_B}_N\subseteq [\rho_A,\rho_B]^{N}$ as the set defined by 
$$
O^{\rho_A,\rho_B}_N:=\left\{\underline m\,:\,\rho_A\leq m_1\leq \dots \leq m_{N}\leq \rho_B\right\}.
$$
The Lebesgue volume is given by $|O^{\rho_A,\rho_B}_N|=\frac{(\rho_B-\rho_A)^{N}}{N!}$.
Our result is the following:
\begin{theorem}\label{ilth}
The invariant measure of the process with generator \eqref{generator} is given by
\begin{equation}\label{formulasuper}
\mu_{N}^{\rho_A,\rho_B}(\eta)=\frac{1}{|O^{\rho_A,\rho_B}_N|}\int_{O^{\rho_A,\rho_B}_N} d \underline m\  \mathcal G_{\underline m}(\eta)\,.
\end{equation}
\end{theorem}
In the above statement we make explicit the dependence of the invariant measure on the parameters $\rho_A, \rho_B, N$, while in the rest of the paper we omit such dependence.
For simplicity of notation we used the symbol $\eta$ for a configuration of particles but
in order to be compatible with our notation for vectors we remark that in \eqref{formulasuper} $\eta\equiv\underline \eta$ should be interpreted as a vector.

In order to better illustrate the result we first give the proof for the case of only one site (Section~\ref{one_site})
and then generalize it for the case of $N$ sites in Section~\ref{N_sites}. The basic telescoping mechanism is active already in the $N=1$ case.

\section{Proof of Theorem \ref{ilth}} \label{sec3}
We introduce the moment generating function of the geometric distribution $\mathcal G_m$:
\begin{equation}\label{gen_fcn}
\mathcal F_m(\lambda):=\sum_{k=0}^{\infty}\mathcal G_m(k)\lambda^k=\left[1+(1-\lambda)m\right]^{-1}\,, \qquad\qquad 0 \le \lambda < \frac{1+m}{m} \,.
\end{equation}
Like before, given $\underline m$ and $\underline \lambda$, we define $\mathcal F_{\underline m}(\underline \lambda):=\prod_{x=1}^N \mathcal F_{m_x}(\lambda_x)$.

\subsection{The case $N=1$}\label{one_site}
In the case that our lattice is composed by one single node which is in contact with two external reservoirs, the state space of the process $\Omega_1$ is the set of natural numbers. We denote by $\eta_1\in \mathbb N_0$ a generic element of the state space and the generator $L_1$ (from \eqref{generator} for $N=1$) is given by
\begin{equation}\label{gen1}
L_1 f(\eta_1)=\sum_{k=1}^{\eta_1}\frac 2k[f(\eta_1-k)-f(\eta_1)]+\sum_{k=1}^{\infty}\frac{\beta_A^k+\beta_B^k}{k}[f(\eta_1+k)-f(\eta_1)]\,,
\end{equation}
where $0<\beta_A <\beta_B <1$ are the parameters associated to the two external reservoirs. 
The stationarity condition for the invariant measure $\mu$ is
\begin{align}\label{stat-1}
 \mu(\eta_1)\left[\sum_{k=1}^{\infty}\frac{\beta_A^k+\beta_B^k}{k}+\sum_{k=1}^{\eta_1}\frac 2k\right]
=\sum_{k=1}^{\infty}\mu(\eta_1+k)\frac 2k+\sum_{k=1}^{\eta_1}\mu(\eta_1-k)\frac{\beta_A^k+\beta_B^k}{k},\nonumber\\
\end{align}
which must be satisfied for all $\eta_1\in\mathbb N_0$.
Theorem \ref{ilth} says that for $N=1$ the invariant measure is a mixture of geometric distributions, i.e.,
\begin{equation}\label{laformula}
\mu(\eta_1)=\frac{1}{\rho_B - \rho_A}\int_{\rho_A}^{\rho_B}dm\ \mathcal G_m(\eta_1)\,.
\end{equation}
Note that in the limit $\rho_A \to \rho_B$ we recover the special equilibrium case,
where the invariant measure is just a geometric distribution of mean $\rho_B$. 

Instead of checking the validity of \eqref{stat-1} for each $\eta_1 \in \mathbb N_0$, it will be convenient to multiply both sides of \eqref{stat-1} by $\lambda^\eta_1$ and sum over $\eta_1$. 
In this way, we get an equality between generating functions for each value of $\lambda$ 
which is equivalent to the whole set of conditions \eqref{stat-1}. 
In the sequel we will use the following elementary formulas:
\begin{align}\label{id2}
 \sum_{k=0}^{\infty}x^k=\frac{1}{1-x}\,; \ \   \sum_{k=1}^{\infty}\frac{x^k}{k}=\log \frac{1}{1-x}\,;\ \ \sum_{k=0}^{+\infty}x^k\sum_{j=1}^k\frac 1j=\frac{1}{1-x}\log\frac{1}{1-x}\, \qquad |x|<1.
\end{align}
We write separately each one of the terms that are obtained by
inserting \eqref{laformula} into \eqref{stat-1} and computing
the generating function. The first term gives 
\begin{align}
\sum_{\eta_1=0}^\infty \lambda^{\eta_1}\mu(\eta_1)\sum_{k=1}^{\infty}\frac{\beta_A^k+\beta_B^k}{k} &
=\frac{1}{\rho_B -\rho_A}\int_{\rho_A}^{\rho_B}\frac{dm}{1+m}\sum_{\eta_1=0}^{\infty}\left(\frac{m\lambda}{1+m}\right)^{\eta_1}\sum_{k=1}^{\infty}\frac{\beta_A^k+\beta_B^k}{k}\\
&=\frac{1}{\rho_B -\rho_A} \int_{\rho_A}^{\rho_B}dm\ \Big[ \log(1+\rho_A) + \log(1+\rho_B)\Big] \mathcal F_m(\lambda)\,.
\end{align}
where we used \eqref{id2}. 
Similarly, exchanging the order of summation, the other terms give
\begin{align}
\sum_{\eta_1=0}^\infty \lambda^{\eta_1} \mu(\eta_1)\sum_{k=1}^{\eta_1}\frac 2k 
&=\frac{1}{\rho_B-\rho_A}\ \int_{\rho_A}^{\rho_B} dm\Big[2\log(1+m)+ 2\log \mathcal F_m(\lambda)\Big] \mathcal F_m(\lambda)\,,
\end{align}
\begin{align}
\sum_{\eta_1=0}^\infty \lambda^{\eta_1}
\sum_{k=1}^{\infty}\mu_1(\eta_1+k)\frac 2k
&=\frac{1}{\rho_B-\rho_A}\ \int_{\rho_A}^{\rho_B} dm\Big[2\log(1+m)\Big] \mathcal F_m(\lambda)
\end{align}
and
\begin{align}
&\sum_{\eta_1=0}^\infty \lambda^{\eta_1}
\sum_{k=1}^{\eta_1}\mu(\eta_1-k)\frac{\beta_A^k+\beta_B^k}{k}\nonumber\\
&=\frac{1}{\rho_B- \rho_A} \int_{\rho_A}^{\rho_B}dm \Big[\log(1+ \rho_A)+\log \mathcal F_{\rho_A}(\lambda)+\log (1+\rho_B)+ \log \mathcal  F_{\rho_B}(\lambda)\Big] \mathcal F_m(\lambda).
\end{align}
All in all, by adding  the terms, we get that the stationarity condition \eqref{stat-1}
is equivalent to  
\begin{equation}\label{semplice}
\int_{\rho_A}^{\rho_B}dm \Big[\log \mathcal F_{\rho_A}(\lambda) -2 \log\mathcal  F_m(\lambda) + \log \mathcal F_{\rho_B}(\lambda)\Big] \mathcal F_m(\lambda) = 0.
\end{equation}
By a direct computation we have the following simple relation for $m\mapsto F_{m}(\lambda)$, the antiderivative of $m\mapsto \mathcal F_{m}(\lambda)$ :
\begin{equation}\label{formulachiave}
{F}_m(\lambda)=\int_0^m dm' \ \mathcal F_{m'}(\lambda)=\frac{1}{(\lambda-1)}\log \mathcal F_m(\lambda)\,.
\end{equation} 
Then, in terms of the antiderivative, \eqref{semplice} is rewritten as
\begin{equation}\label{semplice22}
(\lambda-1) \int_{\rho_A}^{\rho_B}dm \Big[F_{\rho_A}(\lambda) -2  F_m(\lambda) +  F_{\rho_B}(\lambda)\Big]  F'_m(\lambda) = 0,
\end{equation}
where ${F}'_m(\lambda)$ denotes the derivative with respect to the parameter $m$.
Performing the integral, apart the common $(\lambda-1)$ factor, we get 
\begin{equation}
F_{\rho_A}(\lambda) (F_{\rho_B}(\lambda)-F_{\rho_A}(\lambda)) -  (F_{\rho_B}^2(\lambda)-F_{ \rho_A}^2(\lambda)) + F_{\rho_B}(\lambda) (F_{\rho_B}(\lambda)-F_{ \rho_A}(\lambda)),
\end{equation}
which is clearly zero. This concludes the proof of Theorem \ref{ilth} for $N=1$.

\subsection{The general case}\label{N_sites}
In this section we give the proof of Theorem \ref{ilth} for general $N$. 
We now consider the full stationarity condition 
\eqref{stazcondN} which also contains the bulk terms.
With computations similar to the ones done in the previous section we obtain that the stationarity condition \eqref{stazcondN}  is equivalent to 
\begin{align}\label{uno}
\sum_{x=1}^N \int_{O_N^{  \rho_A,\rho_B}} d\underline m \Big[\log \mathcal F_{m_{x-1}}(\lambda_x) -2 \log\mathcal  F_{m_x}(\lambda_x) + \log \mathcal F_{m_{x+1}}(\lambda_x)\Big] \mathcal F_{\underline m}(\underline \lambda) = 0,
\end{align}
where we have defined  
$$\mathcal F_{m_{0}}(\lambda_1) \equiv  \mathcal F_{ \rho_A}(\lambda_1),\qquad \text{and}\qquad
\mathcal F_{m_{N+1}}(\lambda_N) \equiv  \mathcal F_{\rho_B}(\lambda_N).
$$
Using \eqref{formulachiave} the above condition can be also written as
\begin{align}\label{azzerare}
\sum_{x=1}^N (\lambda_x-1) \int_{O_N^{  \rho_A,\rho_B}} d\underline m \Big[ F_{m_{x-1}}(\lambda_x) -2   F_{m_x}(\lambda_x) +  F_{m_{x+1}}(\lambda_x)\Big]  F'_{\underline m}(\underline \lambda) = 0\,,
\end{align}
where, as usual in this paper, we denote $F'_{\underline m}(\underline \lambda)=\prod_{x=1}^NF'_{m_x}(\lambda_x)$.
One can check that the integrals are vanishing for each $x\in\{1,2,\ldots,N\}$. To verify this, 
let us call $O^{  \rho_A,\rho_B}_{N-1,x}$ the collection of $N-1$ ordered variables
$m_y$, with $y\neq x$, i.e., where the variable $m_x$ is missing; we call $\underline{m}^x$ a generic element of $O^{  \rho_A,\rho_B}_{N-1,x}$.
Then, by applying Fubini theorem, we get     
\begin{align}
& \int_{O_N^{  \rho_A,\rho_B}} d\underline m \Big[ F_{m_{x-1}}(\lambda_x) -2   F_{m_x}(\lambda_x) +  F_{m_{x+1}}(\lambda_x)\Big]  F'_{\underline m}(\underline \lambda) \nonumber\\
&=\int_{O_{N-1,x}^{  \rho_A,\rho_B}} d\underline{m}^x\  
\int_{m_{x-1}}^{m_{x+1}}dm_x \Big[ 
{F}_{m_{x-1}}(\lambda_x) F'_{m_x}(\lambda_x)\\
& \hspace{4.cm}-2 F_{m_x}(\lambda_x){F}'_{m_x}(\lambda_x) 
\\
& \hspace{4.cm}
+{F}_{m_{x+1}}(\lambda_x) F'_{m_x}(\lambda_x)
\Big]
F'_{\underline{m}^x}(\underline{\lambda}^x)\,,
\end{align} 
where again $\underline{\lambda}^x$ is obtained from the vector $\underline \lambda$ by removing the component $\lambda_x$.
The integral over the variable $m_x$ on the right hand side of the above equation
can now be performed and we are left with
\begin{align}
\int_{O_{N-1,i}^{  \rho_A,\rho_B}} d\underline{m}^x\  
 \Big[ 
&{F}_{m_{x-1}}(\lambda_x) (F_{m_{x+1}}(\lambda_x) - F_{m_{x-1}}(\lambda_x)) \noindent\\
&-  (F^2_{m_{x+1}}(\lambda_x) - F^2_{m_{x-1}}(\lambda_x))\noindent\\
&+{F}_{m_{x+1}}(\lambda_x) (F_{m_{x+1}}(\lambda_x) - F_{m_{x-1}}(\lambda_x))
\Big]
F'_{\underline{m}^x}(\underline{\lambda}^x)
\end{align} 
which is clearly zero since the term inside the squared parenthesis is identically zero.
This concludes the proof of Theorem \ref{ilth}.

\section{Integrable heat conduction model with parameter $s=1/2$} \label{sec4}
\subsection{The model}\label{ctsmodel}
In this section we show that the same approach based on a direct computation of the joint generating function holds for a related model. The model was introduced  in \cite{FGK} as a scaling limit of the harmonic model and further generalized in \cite{FFG}. The setting is as in the previous section, namely a one dimensional lattice $\Lambda_N$ with two extra ghost sites representing the reservoirs. Here we denote by $z_x\in \mathbb{R}_{+}$ the arbitrary quantity of energy at site $x\in \Lambda_N$, and by $z=(z_1,\dots ,z_N)$ a generic configuration in $\Omega_N = \mathbb{R}_{+}^{\Lambda_N}$, i.e., the state space.
The generator of the stochastic dynamics is given as the superposition of a bulk part and a boundary part, described below:
\begin{equation}\label{generator-levy}
L_N:=L_N^{\textrm{bulk}}+L_N^{\textrm{bound}},
\end{equation}
whose action on functions $f: \Omega_N \to \mathbb{R}$ that are bounded and Lipschitz
 is 
\begin{equation}
\label{genbulk-levy}
L_N^{\textrm{bulk}}f(z)= \sum_{\substack{x,y\in \Lambda_N \\ |x-y|=1}}  \int_{0}^{z_x} \frac{d\alpha}{\alpha} \left[ f\left(z - \alpha \delta_x + \alpha \delta_y \right)  - f\left( z\right)  \right] 
\end{equation} 
and
\begin{align}
\label{genbound-levy}
L_N^{\textrm{bound}}f(z) & = \int_{0}^{z_1} \frac{d\alpha}{\alpha} \left[ f\left(z - \alpha \delta_1  \right)  - f\left( z\right)  \right] +   \int_{0}^{\infty} \frac{d\alpha}{\alpha} e^{-\alpha/T_{A}} \left[ f\left(z + \alpha \delta_1  \right)  - f\left( z\right)  \right]\\  & +
 \int_{0}^{z_N} \frac{d\alpha}{\alpha} \left[ f\left(z - \alpha \delta_N  \right)  - f\left( z\right)  \right] +   \int_{0}^{\infty} \frac{d\alpha}{\alpha} e^{-\alpha/T_{B}} \left[ f\left(z + \alpha \delta_N  \right)  - f\left( z\right)  \right] \;,
\end{align}
where we recall that, as in the discrete case, for $x\in \Lambda_N$,  $\delta_x$ is the configuration with $\delta_x(y)=0$ for $y\neq x$ and $\delta_x(x)=1$.
Above $T_{A}$ (respectively, $T_{B}$) is the temperature associated to the left (respectively, right) reservoir whose purpose is to destroy the conservation of energy by imposing heat conduction from one side of the chain to
the other.
When $T_{A} = T_{B} = T$ there is no transport of energy, the model is reversible and its invariant measure is of product type with each marginal being exponential with mean equal to the temperature $T$ of the  external reservoirs. Note that since $1/\alpha$ is not integrable at zero, this is a jump process with a dense set of jumps. We do not address here the delicate issues related to the definition of the process.


\subsection{Invariant measure}
The stationarity condition imposes that the density $\mu$ of the invariant measure satistisfies
 \begin{align}\label{stazcondN-levy}
&0=\int_0^{z_1}\frac{d\alpha}{\alpha}e^{-\frac{\alpha}{T_A}}\left[\mu(z-\alpha \delta_1)-\mu(z)\right]+\mu(z)\int_{z_1}^{+\infty}\frac {d\alpha}{\alpha}e^{-\frac{\alpha}{T_A}} \\
&+\int_0^{z_1}\frac{d\alpha}{\alpha}\left[\mu(z+\alpha \delta_1)-\mu(z)\right]+\int_{z_1}^{+\infty}\frac {d\alpha}{\alpha} \mu(z+\alpha \delta_1) \\
&+\int_0^{z_N}\frac{d\alpha}{\alpha}e^{-\frac{\alpha}{T_B}}\left[\mu(z-\alpha \delta_N)-\mu(z)\right]+\mu(z)\int_{z_N}^{+\infty}\frac{d\alpha}{\alpha} e^{-\frac{\alpha}{T_B}}\\
&+\int_0^{z_N}\frac{d\alpha}{\alpha}\left[\mu(z+\alpha \delta_N)-\mu(z)\right]+\int_{z_N}^{+\infty}\frac {d\alpha}{\alpha} \mu(z+\alpha \delta_N) \\
&+\sum_{\substack{x,y\in \Lambda_N \\ |x-y|=1}}\left\{ \int_0^{z_x}\frac{d\alpha}{\alpha}\left[\mu(z+\alpha \delta_x-\alpha\delta_y)-\mu(z)\right]+\int_{z_x}^{+\infty}\frac {d\alpha}{\alpha} \mu(z+\alpha \delta_x-\alpha \delta_y)\right\}\,.
\end{align}
\noindent
Let $\mathcal E_m (z)=\frac{1}{m} e^{- z / m} \one_{\{ z \geq 0 \}}$ be the density of an exponential distribution of mean $m>0$. Given $\underline m=( m_1,\dots , m_{N})$ and $\underline z =(z_1, \dots ,z_{N})$ we denote by
$\mathcal E_{\underline m}(\underline z):=\prod_{x=1}^{N}\mathcal E_{m_x}(z_x)$.
As before we  introduce
$O^{T_A,T_B}_N\subseteq [T_A,T_B]^{N}$ as the set defined by 
$$
O^{T_A,T_B}_N:=\left\{\underline m\,:\,T_A\leq m_1\leq \dots \leq m_{N}\leq T_B\right\}.
$$ 
Our result is the following:
\begin{theorem}\label{ilth2-levy}
The invariant measure of the process with generator \eqref{generator-levy} is given by
\begin{equation}\label{formulacont}
\mu_{N}^{T_A,T_B}(z)=\frac{1}{|O^{T_A,T_B}_N|}\int_{O^{T_A,T_B}_N} d \underline m\  \mathcal E_{\underline m}(z)\,.
\end{equation}
\end{theorem}
Here, again, for simplicity of notation we call $z$ a configuration of energies but 
in order to be compatible with our vector-notation we remark that in \eqref{formulacont} $z\equiv\underline z$ should be interpreted as a vector.
\noindent
The strategy of the proof is similar to the previous one, namely we consider $N=1$ first and then we show the result for a general finite chain of $N$ sites. Below, in order to alleviate the notation for the invariant measure, we drop the dependence on the parameters $T_A$, $T_B$ and $N$.

\section{Proof of Theorem \ref{ilth2-levy}} \label{sec5}
We introduce the moment generating function of the exponential distribution $\mathcal E_m$:
\begin{equation}\label{gen_fcn-levy}
\mathcal F_{m} (t):=\int_{0}^{\infty}dz \ \mathcal E_m(z) e^{t z} = \dfrac{1}{1-t m}\,, \qquad  t < \frac{1}{m} \,
\end{equation}
and we define $\mathcal F_{\underline m}(\underline t):=\prod_{x=1}^N \mathcal F_{m_x}(t_x)$.
\subsection{The case $N=1$}
If the lattice consists of only one site then the Markov generator simplifies as
\begin{eqnarray}\label{gen1-levy}
L_1 f(z_1)
&= &
2 \int_{0}^{z_1}  \frac{d\alpha}{\alpha} \left[ f\left(z - \alpha \delta_1  \right)  - f\left( z \right)  \right] \nonumber \\
&+&   \int_{0}^{\infty} \frac{d\alpha}{\alpha} \left( e^{-\alpha/T_{A}} + e^{-\alpha/T_{B}} \right)  \left[ f\left(z + \alpha \delta_1  \right)  - f\left( z\right)  \right].
\end{eqnarray}
The stationary condition for the  invariant measure $\mu$ reads
\begin{align}\label{stat-1-levy-new}
&\int_{0}^{z_1} \frac{d\alpha}{\alpha} \left( e^{-\alpha/T_{A}} + e^{-\alpha/T_{B}} \right)
(\mu(z_1)-\mu(z_1 - \alpha)) +\int_{z_1}^{\infty} \frac{d\alpha}{\alpha} \left( e^{-\alpha/T_{A}} + e^{-\alpha/T_{B}} \right)
\mu(z_1)
\\ &
=  2 \int_{0}^{z_1}  \frac{d\alpha}{\alpha} (\mu(z_1 + \alpha)-\mu(z_1) )
+2 \int_{z_1}^{\infty}  \frac{d\alpha}{\alpha} \mu(z_1 + \alpha)\,,
\end{align}
which must be satisfied for all $z_1\in\mathbb{R}_{+}$. 
Multiplying both sides by $ e^{tz_1} $, using the representation \eqref{formulacont} for $N=1$ and taking the integral in $dz_1$, we get, as in the previous case, four different terms which can be compactly written as

\begin{align} \label{allinall}
 \int_{T_A}^{T_B} dm  \int_{0}^{\infty} d\alpha \ \frac{ e^{\alpha t } -1 }{\alpha}  \left[  \left( e^{-\alpha / T_A } -  2 e^{-\alpha / m}  + e^{-\alpha / T_B}  \right)     \right]  \mathcal F_{m}(t) = 0 \;.
\end{align}
The inner integrals can be computed using ``Feynman's trick''  which, for $a,b > 0$, leads to
\begin{equation}
\int_{0}^{\infty} \dfrac{e^{-ax} - e^{-bx}}{x}  dx = \log\left( \dfrac{b}{a} \right)  \;,
\end{equation} 
so that we have
\begin{align} \label{afterfeynman}
 \int_{T_A}^{T_B} dm   \left[  \log \mathcal F_{T_A}(t)  -2 \log \mathcal F_{m}(t)   +  \log  \mathcal F_{T_B}(t)\right]  \mathcal F_{m}(t) = 0 .
\end{align}
The key observation regarding $F_m$, the antiderivative of $ m \mapsto \mathcal F_{m}(t)$, is the following
\begin{equation}\label{formulachiavelevy}
{F}_{m}(t)=\int_0^m dm' \ \mathcal F_{m'}(t)= - \dfrac{1}{t} \log \left( 1-t m \right)  =  \dfrac{1}{t} \log \mathcal F_{m} (t) \,.
\end{equation}
This allows to write \eqref{afterfeynman} as
\begin{align}
 \int_{T_A}^{T_B} dm   \left[ {F}_{T_A}(t) -2 {F}_{m}(t)  + {F}_{T_B}(t)  \right] {F}'_{m}(t) = 0, 
\end{align}
where  ${F}'_m (t)$ denotes the derivative with respect to the parameter $m$.
As before, by inspection the left hand side of the previous equation is zero and the proof of Theorem \ref{ilth2-levy} for $N=1$ is concluded.

\subsection{The general case}
For general  $N$ the  stationarity condition is written in equation \eqref{stazcondN-levy}.
As before, we multiply both sides by $\displaystyle \prod_{x=1}^{N} e^{t_x z_x} $, we use the representation \eqref{formulacont} and take the integral. We obtain
\begin{align}
&\int_{O_N^{T_{A},T_{B}}} d\underline m \, \Big[ \int_{0}^{\infty}   \frac{d\alpha}{\alpha} \left(  e^{- \alpha/T_A}  + \sum_{x=1}^{N}  2   e^{\alpha t_x} e^{-\alpha/m_x} + e^{-\alpha/T_B}  \right) \mathcal F_{\underline m}(\underline t)  \Big] = \nonumber\\ & 
\label{hi}
\int_{O_N^{T_{A},T_{B}}}  d\underline m \,  \int_{0}^{\infty}   \frac{d\alpha}{\alpha} \left[ e^{- \alpha / m_1}    + \sum_{x=1}^{N}   e^{\alpha t_x} \left( e^{-\alpha / m_{x-1}} + e^{-\alpha / m_{x+1}} \right)   + e^{- \alpha / m_N}  \right] \mathcal F_{\underline m}(\underline t),
\end{align}
where we have set $m_0 := T_A$ and  $m_{N+1} := T_B$.
At this point it is enough to notice that using the telescoping cancellation
\begin{align*}
 \int_{0}^{\infty}   \frac{d\alpha}{\alpha} & \left(  e^{- \alpha/T_A}  -   e^{- \alpha/m_1} - e^{- \alpha / m_N} + e^{- \alpha/T_B}  \right)  \\  & = \sum_{x=1}^{N} \int_{0}^{\infty}   \frac{d\alpha}{\alpha} \left(  e^{- \alpha/m_{x-1}} - 2 e^{- \alpha/m_x} + e^{- \alpha/m_{x+1}}  \right) 
\end{align*}
we can rewrite \eqref{hi} in a form analogous to  \eqref{allinall}:
\begin{equation}
\sum_{x=1}^{N} \int_{O_N^{T_A,T_B}} d\underline m \,    \int_{0}^{\infty}   {d\alpha}\left( \dfrac{e^{\alpha t_x} -1}{\alpha} \right)  \left[ e^{- \alpha/ m_{x-1}}-2 e^{- \alpha/ m_{x}}+ e^{- \alpha/ m_{x+1}} \right]  \mathcal F_{\underline m}(\underline t) = 0 \;.
\end{equation}
Computing the inner integrals, we get
\begin{align}
\sum_{x=1}^{N} \int_{O_N^{T_A,T_B}} d\underline m \,     \left[  \log \mathcal F_{m_{x-1}}(t_x)-2 \log \mathcal F_{m_x}(t_x)  +  \log  \mathcal F_{m_{x+1}}(t_x)   \right]  \mathcal F_{\underline m}(\underline t)  = 0 ,
\end{align}
which can be written in terms of the antiderivative $F_{\underline m}$ using equation \eqref{formulachiavelevy}
\begin{align}
\sum_{x=1}^{N} \int_{O_N^{T_A,T_B}} d\underline m \,     \left[ F_{m_{x-1}} (t_x) -2 F_{m_x}(t_x)   + F_{m_{x+1}} (t_x) \right]   F'_{\underline m}(\underline t)  = 0 \;.
\end{align}
We show that each term of the above sum is zero.  To this aim we apply Fubini theorem to the $x^{th}$ term to separate the integral in $m_x$, i.e.,
\begin{align*}
& \int_{O_N^{T_A,T_B}} d\underline m \,     \left[ F_{m_{x-1}} (t_x) -2 F_{m_x}(t_x) + F_{m_{x+1}} (t_x)   \right]   F'_{\underline m}(\underline t)  =  \\  &
\int_{O_{N-1,x}^{T_A,T_B}} d\underline m^{x} \,  \int_{m_{x-1}}^{m_{x+1}} dm_x \left[ F_{m_{x-1}} (t_x) -2 F_{m_x}(t_x) + F_{m_{x+1}} (t_x)  \right] F'_{m_x} (t_x) \prod_{ \substack{y=1 \\  y \neq x}}^{N}F'_{m_y} (t_y),
\end{align*}
where $O_{N-1,x}^{T_A,T_B}$ has the same meaning as before, namely the collection of $N-1$ ordered variables $m_y$ with $y \neq x$. 
Computing the inner integral on the right hand side we obtain zero and the proof of Theorem \ref{ilth2-levy} is concluded.

\section{Some applications}
In this last section we discuss some important applications
that follow from the representation of the invariant measure as a mixture. 
%

First, we deduce a general FKG-type inequality.
In particular, we show that when sampled according to the steady states of Theorem \ref{ilth} and Theorem \ref{ilth2-levy}, the processes are {\it associated} in the sense of Definition 1.1 \cite{Esary}, which we recall below. 
To this end we define a partial ordering in $\mathbb R_+^{\Lambda_N}$ by saying that $\underline X \le \underline Y$ if for all $i\in \Lambda_N$,  $X_i \le Y_i$. Then a function $g: \mathbb R^{\Lambda_N} \to \mathbb R$ is said {\it non-decreasing} if, for all pairs $\underline X, \underline Y\in \mathbb R^{\Lambda_N}$ with $\underline X \le \underline Y$, we have $g(\underline X)\le g(\underline Y)$. A random variable $\underline X \in \mathbb R_+^{\Lambda_N}$ is called associated if for all non-decreasing functions $g,h: \mathbb R_+^{\Lambda_N} \to \mathbb R$, 
\be
\mathbb E[g(\underline X) h(\underline X)] \ge \mathbb E[g(\underline X)]\cdot \mathbb E[h(\underline X)] \,
\ee
and the same terminology is used for the corresponding distribution of the random variable $\underline X$.

In the following we show that, thanks to the representation as a mixture of the stationary measures of the processes, we deduce in few steps the association of those measures relying on the well establised association of the ordered statistics of i.i.d random variables given in \cite{Esary}.
The statement of the following theorem was suggested by an anonymous referee.
\begin{theorem}\label{ref}
The invariant measure  $\mu_N^{\rho_A, \rho_B}$ of the model defined in Section \ref{sec11} is associated, namely
for all non-decreasing functions $
g,h:\Omega_N \to \mathbb R$,  we have
\be
\mathbb{E}_{\mu_N^{\rho_A, \rho_B}}  \left[ g(\eta) h(\eta)  \right] \ge \mathbb{E}_{\mu_N^{\rho_A, \rho_B}} \left[  g(\eta) \right]  \cdot \mathbb{E}_{\mu_N^{\rho_A, \rho_B}} \left[   h(\eta) \right].
\ee
The same property holds for the invariant measure $\mu_N^{T_A, T_B}$ of the process defined in Section \ref{ctsmodel}.
\end{theorem}
\begin{proof}
We have to prove that 
\begin{align}
&\frac{1}{|O^{\rho_A,\rho_B}_N|} \int_{O^{\rho_A,\rho_B}_N} d \underline m \sum_\eta g(\eta)h(\eta)\mathcal G_{\underline m}(\eta)\label{one}  \\
&\ge \frac{1}{|O^{\rho_A,\rho_B}_N|} \bigg( \int_{O^{\rho_A,\rho_B}_N} d \underline m \sum_\eta g(\eta) \mathcal G_{\underline m}(\eta)\bigg) \cdot \frac{1}{|O^{\rho_A,\rho_B}_N|} \bigg( \int_{O^{\rho_A,\rho_B}_N} d \underline m \sum_\eta h(\eta)\mathcal G_{\underline m}(\eta) \bigg).
\end{align}
First of all we observe that, from Theorem 2.1 in \cite{Esary}, we have that any product measure is associated
thus, for all fixed $\underline m$, 
\begin{align}
 \sum_\eta g(\eta)h(\eta)\mathcal G_{\underline m}(\eta) 
\ge \bigg(  \sum_\eta g(\eta) \mathcal G_{\underline m}(\eta)\bigg) \cdot \bigg( \sum_\eta h(\eta) \mathcal G_{\underline m}(\eta) \bigg) \,.
\end{align}
We define now the functions $\tilde g, \tilde h: \mathbb R_+^{\Lambda_N}\to \mathbb R$ as
\be
\tilde g(\underline m):=\sum_\eta g(\eta)\mathcal G_{\underline m}(\eta)
\ee
and
\be
\tilde h(\underline m):=\sum_\eta h(\eta) \mathcal G_{\underline m}(\eta) \ .
\ee
Since $\underline m \le \underline m'$ implies  $\mathcal G_{\underline m}\preccurlyeq \mathcal G_{\underline m'}$ (where the symbol $\preccurlyeq$ indicates stochastic domination), it follows that  $\tilde g$ and $\tilde h$ are non-decreasing functions. 
As a consequence, in order to prove \eqref{one} it is sufficient to show that
\begin{align}
 \int_{O^{\rho_A,\rho_B}_N}  \tilde g(\underline m) \cdot \tilde h (\underline m) d \underline m \ge \frac{1}{|O^{\rho_A,\rho_B}_N|} \bigg( \int_{O^{\rho_A,\rho_B}_N}  \tilde g(\underline m) d \underline m \bigg) \cdot \bigg( \int_{O^{\rho_A,\rho_B}_N} \tilde h(\underline m) d \underline m \bigg) 
\end{align}
which follows from the association of ordered statistics of i.i.d. random variables, as shown in Section 5 of \cite{Esary}.
Notice that also for an exponential distribution of mean $m$, it holds that if $\underline m \le \underline m'$ then   $\mathcal E_{\underline m}\preccurlyeq \mathcal E_{\underline m'}$ and so the statement is also true for the continuous model.
\end{proof}

We remark that positive correlation inequalities have also been obtained \cite{GKF} for other models having convex quadratic mobility,
such as the symmetric inclusion process and the Brownian energy process.
It could be interesting to further investigate whether association is true as well.

Another consequence of the representation of the invariant measure as a mixture of independent random variables is the proof of the large deviation principle for the density profile
which can be deduced by a combination of two large deviation principles: one for the order statistics and another for independent inhomogeneous random variables with an additional application of the contraction principle. The heuristic argument is outlined in \cite{BDGJL-int}, Section 3.2.
A rigorous proof is given
in \cite{CFFGR}, where it is computed the pressure, the density large deviation functional 
and their additivity principle. This provided a rigorous proof for the expression of the density large deviation function
for the whole class of harmonic models obtaining a rate function in accordance 
with the result
of the MFT \cite{MFT} for systems with convex quadratic mobility and constant diffusion.



\section*{Declarations}
The authors have no conflicts of interest. Data sharing is not applicable to this article as no datasets were generated or analysed during the current study.

\vskip.2cm
\noindent
{\bf Acknowledgments.} We thank an anonymous referee for suggesting the statement of Theorem \ref{ref}.
GC, CF and CG acknowledge several useful discussions with R. Frassek and F. Redig about the
non-equilibrium steady state of harmonic models.
We acknowledge financial support from the Italian Research Funding Agency (MIUR) through
PRIN project ``Emergence of condensation-like phenomena in interacting particle systems: kinetic and lattice models'', grant n. 202277WX43. We also thank Istituto Nazionale di Alta Matematica (INDAM) for its support.

\end{document}